\documentclass{svmult}

\usepackage{makeidx}         
\usepackage{graphicx}        
\usepackage{multicol}        
\usepackage[bottom]{footmisc}

\makeindex             

\newtheorem{thm}{Theorem}[section]
\newtheorem{lem}{Lemma}[section]
\newtheorem{prop}{Proposition}[section]
\newtheorem{cor}{Corollary}[section]
\newtheorem{rmk}{Remark}[section]
\newtheorem{exap}{Example}[section]

\begin{document}

\title*{Generalized It$\hat {\rm o}$ Formulae and Space-Time Lebesgue-Stieltjes Integrals of Local Times: Corrected\thanks{Originally published in {\it  S\'eminaire de
Probabilit\'es}, Vol.40 (2007), 117-136. There was an error in the statement of the conditions of Theorem 2.1, which we had hoped to correct ``in proof",  
leading to this unpublished version (of 1 March 2007). The following two papers have further developments 
\newline
{[1]} C.R. Feng and H.Z. Zhao, Rough path integral of local time, {\it Comptes rendus de l'Academie des sciences Paris, Ser. I}, Vol. 346 (2008), 431-434. 
\newline
{[2]} C.R. Feng and Z.Z. Zhao, Local time rough path for  L\'evy processes,
{\it Electronic Journal of Probability}, Vol. 15 (2010), 452-483. 
}}
\titlerunning{Generalized It$\hat {\rm o}$ Formulae }
\author{K.D. Elworthy\inst{1},
A. Truman\inst{2} and H.Z. Zhao\inst{3}}
\institute{$^1$ Mathematics Institute, University of Warwick, Coventry CV4 7AL,UK.
\texttt{kde@maths.warwick.ac.uk}
\and $^2$ Department of Mathematics,
University of Wales Swansea,
Singleton Park, Swansea SA2 8PP, UK.
\texttt{A.Truman@swansea.ac.uk}
\and $^3$ Department of Mathematical Sciences, Loughborough
University, LE11 3TU, UK.
\texttt{H.Zhao@lboro.ac.uk}}
%
%
\maketitle

\renewcommand{\theequation}{\arabic{section}.\arabic{equation}}

\newcounter{bean}

\begin{abstract}\
Generalised It${\hat{\rm o}}$ formulae are proved for time dependent
functions of continuous real valued semi-martingales.The conditions involve
left space and time first  derivatives, with the left space derivative required to have locally
bounded 2-dimensional variation. In particular a class of functions with discontinuous first
derivative is included. An estimate of Krylov allows further weakening of these conditions
when the semi-martingale is a diffusion.
\vskip5pt

\vskip5pt


\end{abstract}

\section{Introduction}
Extensions of It${\hat {\rm o}}$ formula to less smooth functions are useful in
studying many problems such as partial differential equations with some
singularities, see below, and in the mathematics of finance.
The first extension was obtained for $|X(t)|$ by Tanaka \cite{tanaka} with a beautiful use of local time.
The generalized It${\hat {\rm o}}$ formula in one-dimension for time independent
convex functions was developed in \cite{meyer} and for superharmonic functions in multidimensions
in \cite{brosamler} and for distance functions in \cite{kendall}. Extensions of  It${\hat {\rm o}}$'s
formula have also been studied by \cite{krylov}, \cite{follmer}, \cite{nualart} and \cite{frw}.
In \cite{frw},  It${\hat {\rm o}}$'s formula for $W^{1,2}_{loc}$ functions was
studied using Lyons-Zheng's backward and forward stochastic integrals
\cite{lz}. In \cite{yor1}, It${\hat {\rm o}}$'s formula was extended to absolutely continuous
functions with locally bounded derivative using the integral $\int _{-\infty}^{\infty}
\nabla f(x){\rm d}_xL  _s(x)$. This integral was defined
through the existence of the expression  $f(X(t))-f(X(0))-\int _0^t{\partial \over
\partial x}f(X(s)){\rm d}X(s)$; it was extended to $\int_0^t\int _{-\infty}^{\infty}
\nabla f(s,x){\rm d}_{s,x}L  _s(x)$ for a time dependent function $f(s,x)$ using forward and backward
integrals  for Brownian motion in \cite{eisenbaum}.  Recent activities in this
direction have been to look for minimal assumptions on $f$ to make
this integral well defined for semi-martingales other
than Brownian motion \cite{eisenbaum2}. However, our motivation in
establishing generalized It$\hat {\rm o}$ formulae was to use them to describe the
 asymptotics of the solution of heat equations in the presence of a caustic. Due to the
 appearance of caustics, the solution of the Hamilton-Jacobi equation, the leading term in the asymptotics, is no longer differentiable, but has a jump in the gradient  across
 the shock wave front of the associated Burgers' equation.  Therefore,
the local time of continuous semimartingales in a neighbourhood of the shock wave front of the Burgers equation and the jump of the derivatives of the Hamilton-Jacobi function
(or equivalently the jump in the Burgers' velocity) appear naturally in the semi-classical representation of the corresponding solution to the heat equation \cite{etz}. None of the earlier versions of
 It$\hat {\rm o}$'s formula apply directly to this situation.

In this paper, we first  generalize It${\hat {\rm o}}$'s formula to the case of a continuous semimartingale and a function $f(t,x)$ which is absolutely continuous in each variable and satisfies (1) its left derivative ${\partial ^- \over \partial t}f(t,x)$ exists and
 is  left continuous, (2)
$f(t,x)=f_h(t,x)+f_v(t,x)$ with $f_h(t,x)$ being $C^1$ in $x$ and
 $\nabla f_h(t,x)$ having left continuous and locally bounded left derivative
 $\Delta ^-f_h(t,x)$, and $f_v$ having
left derivative $\nabla  ^-f_v(t,x)$ which is left continuous and of locally bounded variation in $(t, x)$.
Here we use the two-dimensional Lebesgue-Stieltjes
integral of local time with respect to $\nabla^- f(t,x)$.
The main result of this paper is formula (\ref{krylov6}). Formula
(\ref{ito11}) follows from (\ref{krylov6}) easily as a special case.
These formulae appear to be new and in a good form for extensions to
two dimensions (\cite{fz1}). Moreover, in \cite{fz2}, Feng and Zhao observed that the local time
$L_t(x)$ can be considered as a  rough path in $x$ of finite
 2-variation and therefore defined $\int _0^t\int_{-\infty}
^{\infty} \nabla ^-f(s,x){\rm {\rm d}_{s,x}}L_s(x)$ pathwise by extending Young and Lyons' profound
 idea of rough path
integration (\cite{lq}, \cite{young}) to two parameters.
When this paper was nearly completed, we received two preprints concerning a generalized
It${\hat {\rm o}}$'s formula for a continuous function $f(t,x)$ with jump
derivative $\nabla^- f(t,x)$, (\cite{pe}, \cite{gp}). We remark that
 formula (\ref{ito11}) was also observed by \cite{pe} independently.

 In section 3, we consider diffusion processes $X(t)$. We prove the generalized It${\hat{\rm o}}$
 formula for a function $f$ with
generalized derivative ${\partial\over \partial t} f$ in $L^2_{loc}(dtdx)$ and generalized derivative $\nabla f(t,x)$ being of locally
bounded variation in $(t,x)$. We use an inequality from Krylov \cite{krylov}.

\section{The continuous semimartingale case}

We need the following definitions (see e.g. \cite{ash}, \cite{mcshane}): A two-variable function $f(s,x)$ is called monotonically increasing if whenever $s_2\geq s_1$, $x_2\geq x_1$,
\begin{eqnarray*}
f(s_2,x_2)-f(s_2,x_1)-f(s_1,x_2)+f(s_1,x_1)\geq 0.
\end{eqnarray*}
It is called monotonically decreasing if $-f$ is monotonically increasing. The function $f$ is called left continuous iff it is left continuous in both variables together, in other words, for any sequence
$(s_1,x_1)\leq (s_2,x_2)\leq \cdots \leq (s_k,x_k)\to (s,x)$, we have $f(s_k,x_k)\to f(s,x)$ as
$k\to \infty$. Here $(s,x)\leq (t,y)$ means $s\leq t$ and $x\leq y$.
For a monotonically
increasing and left continuous function $f(s,x)$, we can define a Lebesgue-Stieltjes measure by setting
\begin{eqnarray*}
\mu ([s_1,s_2)\times [x_1,x_2))=f(s_2,x_2)-f(s_2,x_1)-f(s_1,x_2)+f(s_1,x_1),
\end{eqnarray*}
for $s_2>s_1$ and $x_2>x_1$. So for a measurable function $g(s,x)$, we can define the
Lebesgue-Stieltjes integral by
\begin{eqnarray*}
\int _{t_1}^{t_2}\int _a^b g(s,x){\rm d}_{s,x} f(s,x)=\int _{t_1}^{t_2}\int _a^b g(s,x)d\mu.
\end{eqnarray*}
Denote a partition ${\cal P}$ of $[t,s]\times
[a,x]$ by $t=s_1<s_2<\cdots <s_m=s$, $a=x_1<x_2<\cdots <x_n=x$ and the variation of $f$ associated with ${\cal P}$  by
\begin{eqnarray*}
&&
V_{\cal P}(f,[t,s]\times [a,x])\\
&=&\sum _{i=1}^{m-1}\sum _{j=1}^{n-1} |f(s_{i+1},x_{j+1})-f(s_{i+1},x_j)-f(s_i,x_{j+1})+f(s_i,x_j)|
\end{eqnarray*}
and the variation of $f$ on $[t,s]\times [a,x]$ by
\begin{eqnarray*}
V_f([t,s]\times [a,x])=\sup _{\cal P}V_{\cal P}(f,[t,s]\times [a,x]).
\end{eqnarray*}

One can find Proposition \ref{prop2}, its proof and definition of the multidimensional
Lebesgue-Stieltjes integral with respect to measures generated by functions of bounded variation
in \cite{mcshane}.
For the convenience of the reader, we include them here briefly.
\begin{prop} \label {prop1}(Additivity of variation) For $s_2\geq s_1\geq t$, and $a_2\geq a_1\geq a$,
\begin{eqnarray}
V_f([t,s_2]\times [a,a_2])&=&V_f([t,s_1]\times [a,a_2])+V_f([t,s_2]\times [a,a_1])\nonumber\\
&&+V_f([s_1,s_2]\times [a_1,a_2])-V_f([t,s_1]\times [a,a_1]).
\end{eqnarray}
\end{prop}
\vskip5pt
{\em Proof}. We only need to prove that for $a\leq a_1< a_2$ and $t\leq s_1$,
\begin{eqnarray}\label{variation1}
V_f([t,s_1]\times [a,a_2])=V_f([t,s_1]\times [a,a_1])+V_f([t,s_1]\times [a_1,a_2]).
\end{eqnarray}
Our proof is similar to the case of one-dimension. We can always refine a partition ${\cal P}$ of
$[t,s_1]\times [a,a_2]$ to include $a_1$. The refined partition is denoted by
${\cal P}^{\prime}$. Then
\begin{eqnarray*}
V_{\cal P}(f,[t,s_1]\times [a,a_2])\leq V_{{\cal P}^{\prime}}(f,[t,s_1]\times [a,a_2]).
\end{eqnarray*}
Then (\ref{variation1}) follows easily. \hfill$\diamond$

\begin{prop}\label{prop2}
A function $f(s,x)$ of locally bounded variation can be decomposed as the
difference of two increasing functions $f_1(s,x)$ and $f_2(s,x)$, in any quarter space
$s\geq t, x\geq a$.
Moreover, if $f$ is also left
continuous, then $f_1$ and $f_2$ can be taken left continuous.
\end{prop}
\vskip5pt
{\em Proof}. For any $(t,x)\in R^2$, define for $s\geq t$ and $x\geq a$,
\begin{eqnarray*}
2\tilde f_1(s,x)&=&V_f([t,s]\times [a,x])+f(s,x),\\
2\tilde f_2(s,x)&=&V_f([t,s]\times [a,x])-f(s,x).
\end{eqnarray*}
Then $f(s,x)=\tilde f_1(s,x)-\tilde f_2(s,x)$. We need to prove that $\tilde f_1$ and $\tilde
f_2$ are increasing functions. For this,
let $s_2\geq s_1\geq t$, $a_2\geq a_1\geq a$,
then use Proposition \ref{prop1},
\begin{eqnarray*}&&
2(\tilde f_1(s_2,a_2)-\tilde f_1(s_1,a_2)-\tilde f_1(s_2,a_1)+\tilde f_1(s_1,a_1))\\
&=&V_f([t,s_2]\times [a,a_2])-V_f([t,s_1]\times [a,a_2])-V_f([t,s_2]\times [a,a_1])\\
&&
+V_f([t,s_1]\times [a,a_1])+f(s_2,a_2)-f(s_1,a_2)-f(s_2,a_1)+f(s_1,a_1)\\
&=&V_f([s_1,s_2]\times [a_1,a_2])+f(s_2,a_2)-f(s_1,a_2)-f(s_2,a_1)+f(s_1,a_1)\\
&\geq & 0.
\end{eqnarray*}
So $\tilde f_1(s,x)$ is an increasing function. Similarly one can prove that $\tilde
f_2(s,x)$ is an increasing function.

Define
\begin{eqnarray*}
f_i(s,x)&=&\lim \limits_{t\uparrow s,y\uparrow x}\tilde f_i(t,y),\ \ i=1,2.
\end{eqnarray*}
Then since $f$ is left continuous, so
\begin{eqnarray}\label{andrew1}
f(s,x)=f_1(s,x)-f_2(s,x),
\end{eqnarray}
and $f_1$ and $f_2$ are as required.
\hfill $\diamond$
\bigskip


From Proposition \ref{prop2},
the two-dimensional Lebesgue-Stieltjes integral of a measurable function $g$
with respect to the left continuous function $f$ of bounded variation
can be defined by
\begin{eqnarray*}
\int _{t_1}^{t_2}\int _a^b g(s,x){\rm d}_{s,x} f(s,x)&=&\int _{t_1}^{t_2}\int _a^b g(s,x){\rm d}_{s,x} f_1(s,x)
\\
&&
-\int _{t_1}^{t_2}\int _a^b g(s,x){\rm d}_{s,x} f_2(s,x)\  {\rm for } \  t_2\geq t_1, b\geq a.
\end{eqnarray*}
Here $f_1$ and $f_2$ are taken to be left continuous.
\smallskip

It is worth pointing out that it is possible that a function $f(s,x)$ is of locally bounded variation in
$(s,x)$ but not of locally bounded variation in $x$ for fixed $s$.
For instance consider $f(s,x)=b(x)$, where $b(x)$ is not of locally
bounded variation, then $V_f=0$. However it is easy to see that when a function
$f(s,x)$ is of locally bounded variation in $(s,x)$ and of locally bounded variation
in $x$ for a fixed $s=s_0$, then it is of locally bounded variation in $x$
for all $s$. We denote by $V_{f(s)}[a,b]$ the variation of $f(s,x)$ on $[a,b]$
as a function of $x$ for a fixed $s$.
\smallskip

Now we recall some well-known results of local time which will be used later in this paper.
Let $X(s)$ be a continuous semimartingale $X(s)=X(0)+M_s+V_s$ on a probability space $\{\Omega,{\cal F},P\}$. Here $M_s$ is a
continuous local martingale and $V_s$ is a continuous process of bounded variation. Let $L _t(a)$ be the local time introduced by P. L\'evy
\begin{eqnarray}
L  _t(a)=\lim_{\epsilon\downarrow 0} {1\over 2\epsilon}\int _0^t1_{[a,a+\epsilon)}(X(s))d<M,M>_s \ \
a.s.,
\end{eqnarray}
for each $t$ and $a$.
 Then it is well known that for each fixed $a\in R$,
 $L _t(a,\omega)$
is continuous, and nondecreasing in $t$ and right continuous with left limit (cadlag) with respect to $a$
(\cite{ks}, \cite{yor}). Therefore we can consider the Lebesgue-Stieltjes integral
$\int _0^{\infty}\phi(s)dL  _s(a,\omega)$
for each $a$ for any Borel-measurable function $\phi$. In particular
\begin{eqnarray}
\int _0^{\infty}1_{R-\{a\}}(X(s))dL _s(a,\omega)=0 \  \ a.s.
\end{eqnarray}
Furthermore if $\phi$ is in $L^{1,1}_{loc}(ds)$, i.e. $\phi$ has locally integrable generalized
derivative,
then we have the following integration by parts formula
\begin{eqnarray}
\int _0^t\phi(s)dL _s(a,\omega)=\phi(t)L _t(a,\omega)-\int _0^t\phi^{\prime}(s)L _s(a,\omega)ds\ \ a.s.
\end{eqnarray}
Moreover, if $g(s,x)$ is Borel measurable in $s$ and $x$ and bounded, by the occupation times formula (e.g. see \cite{ks}, \cite{yor})),
\begin{eqnarray*}
\int _0^tg(s,X(s))d<M,M> _s=2\int _{-\infty}^{\infty}\int _0^tg(s,a)dL _s(a,\omega)da\ \ a.s.
\end{eqnarray*}
If further $g(s,x)$ is in $L^{1,1}_{loc}(ds)$ for almost all $x$, then using the integration by parts formula, we have
\begin{eqnarray*}
\int _0^tg(s,X(s))d<M,M> _s&=&2\int _{-\infty}^{\infty}\int _0^tg(s,a)dL _s(a,\omega)da\nonumber\\
&=&2\int _{-\infty}^{\infty}g(t,a)L _t(a,\omega)da\nonumber\\
&&-2\int _{-\infty}^{\infty}\int _0^t{\partial \over \partial s}g(s,a)L _s(a,\omega)dsda\ \ a.s.
\end{eqnarray*}

We first prove a theorem with $f_h= 0$. The result with a term $f_h$ is a trivial
generalization of Theorem \ref{tom100}.

\begin{thm} \label{tom100}
Assume $f: [0,\infty)\times R\to R$ satisfies
\vskip2pt

(i) $f$ is absolutely continuous in each variable,
\vskip2pt

(ii) the left derivatives ${\partial ^-\over \partial t}f$ and  $\nabla ^-f$ exist at all points of $(0,
\infty )\times R$ and $[0,\infty)\times R$ respectively,
\vskip2pt

(iii)  ${\partial ^-\over \partial t}f$ and  $\nabla ^-f$ are left continuous and locally bounded,
\vskip2pt

(iv)  $\nabla ^-f$ is of locally bounded variation in $(t, x)$ and $\nabla ^-f(0,x)$ is of locally
bounded variation in $x$.
\vskip2pt

 \noindent
 Then for any continuous semimartingale $\{X(t), t\geq 0\}$,
 $f(t,X(t))$ is a semi-martingale and
 \begin{eqnarray}\label{ito10}
&&
 f(t,X(t))\nonumber\\
&=&f(0,X(0))+\int _0^t{\partial^- \over \partial s} f(s,X(s))ds+\int _0^t\nabla ^-f(s,X(s))dX(s)\\
&&
+
\int _{-\infty}^{\infty }L  _t(x){\rm d}_x\nabla ^- f(t,x)-\int _{-\infty}^{+\infty}\int _0^{t}L  _s(x)
{\bf \rm d}_{s,x}\nabla ^-f(s,x)\ \ a.s. \nonumber
\end{eqnarray}
\end{thm}
 \vskip5pt
 {\em Proof}. By a standard localization argument we can assume that $X$ and its quadratic
 variation are bounded processes and that $f$,  ${\partial ^-\over \partial t}f$, $\nabla ^-f$,
 $V_{\nabla ^-f(t)}$ and $V_{\nabla ^-f}$ are bounded (note here $V_{\nabla ^-f(0)}<\infty$ and
 $V_{\nabla ^-f}<\infty$ imply $V_{\nabla ^-f(t)}<\infty$ for all $t\geq 0$).
 We use standard regularizing mollifiers (e.g. see \cite{ks}).
Define
\begin{eqnarray*}
\rho(x)=\cases {c{\rm e}^{{1\over (x-1)^2-1}}, {\rm \  if } \ x\in (0,2),\cr
0, \ \ \ \ \ \ \ \ \ \ \ {\rm otherwise.}}
\end{eqnarray*}
Here $c$ is chosen such that $\int _0^2\rho(x)dx=1$.
Take $\rho_n(x)=n\rho(nx)$ as mollifiers. Define
\begin{eqnarray}\label{krylov5}
f_n(s,x)=\int _{-\infty}^{+\infty}\int _{-\infty}^{+\infty}\rho_n(x-y)\rho_n(s-\tau)f(\tau,y)d\tau dy, \ \ n\geq 1,
\end{eqnarray}
 where we set $f(\tau,y)=f(-\tau,y)$ if $\tau<0$.
 Then $f_n(s,x)$ are smooth and
\begin{eqnarray}\label{truman41}
f_n(s,x)=\int _0^2\int _0^2\rho(\tau)\rho(z)f(s-{\tau\over n},x-{z\over n})d\tau dz, \ \ n\geq 1.
\end{eqnarray}
Because of the absolutely continuity, we can differentiate
under the integral in (\ref{truman41}) to see
that  ${\partial \over \partial t} f_n(t,x)$,
$\nabla f_n(t,x)$, $V_{\nabla f_n(t)}$ and $V_{\nabla f_n}$ are uniformly bounded. In particular
\begin{eqnarray}
\nabla f_n(t,x)=\int _0^2\int _0^2\rho(\tau)\rho(z)\nabla ^-f(s-{\tau\over n},x-{z\over n})d\tau dz
,\ \ t\geq 0.\label{zhao33}
\end{eqnarray}

Moreover
using Lebesgue's dominated convergence theorem, one can prove that as $n\to \infty$,
for each $(t,x)$ with $t\geq 0$,
\begin{eqnarray}
f_n (t,x)&\to & f (t,x).\label{zhao1}
\end{eqnarray}
Also
\begin{eqnarray}
{\partial \over \partial t}f_n(t,x)&\to & {\partial  ^-\over \partial t} f(t,x),\ \ t>0\label{zhao2}\\
\nabla f_n(t,x)&\to&  \nabla ^-f(t,x),\ \ t\geq 0.\label{zhao3}
\end{eqnarray}
Note the convergence
in (\ref{zhao1}), (\ref{zhao2}), (\ref{zhao3}) is also in $L^p_{loc}$, $1\leq p<\infty$.

Now we can use It${\hat {\rm o}}$'s formula for the smooth function $f_n(s,X(s))$, then a.s.
\begin{eqnarray}\label{zhao11}
f_n(t,X(t))-f_n(0,X(0))&=&\int _0^t{\partial \over \partial s} f_n(s,X(s))ds+\int _0^t\nabla f_n(s,X(s))dX(s)\nonumber\\
&&
+{1\over 2}
\int _0^t\Delta f_n(s,X(s))d<M,M>_s.
\end{eqnarray}
As $n\to \infty$, for all $t\geq 0$,
\begin{eqnarray*}
f_n(t,X(t))-f_n(0,X(0))\to f(t,X(t))-f(0,X(0))\ \ a.s.,
\end{eqnarray*}
and
\begin{eqnarray*}
\int _0^t{\partial \over \partial s} f_n(s,X(s))ds\to \int _0^t{\partial ^-\over \partial s} f(s,X(s))ds\ \ a.s.,
\end{eqnarray*}
\begin{eqnarray*}
\int _0^t \nabla f_n(s,X(s))dV_s\to \int _0^t\nabla^- f(s,X(s))dV_s\ \ a.s.
\end{eqnarray*}
and
\begin{eqnarray*}
E\int _0^t(\nabla f_n(s,X(s)))^2d<M,M>_s\to E \int _0^t(\nabla^- f(s,X(s)))^2d<M,M>_s.
\end{eqnarray*}
Therefore in $L^2(\Omega,P)$,
\begin{eqnarray*}
\int _0^t \nabla f_n(s,X(s))dM_s\to \int _0^t\nabla^- f(s,X(s))dM_s.
\end{eqnarray*}
To see the convergence of the last term, we recall the well-known result that the local time
$L  _s(x)$ is jointly continuous in $s$ and cadlag with respect to $x$ and
has a compact support in space $x$ for each $s$.  Now we use the occupation times formula,
\begin{eqnarray}
&&
{1\over 2}\int _0^t\Delta f_n(s,X(s))d<M,M>_s\nonumber\\
 &=&\int _{-\infty}^{+\infty}\int _0^t\Delta f_n(s,x){\rm d}_sL  _s(x)dx\nonumber\\
&=&
\int _{-\infty}^{+\infty}\Delta f_n(t,x)L  _t(x)dx-\int _{-\infty}^{+\infty}\int _0^{t}{d \over ds}
\Delta f_n(s,x)L  _s(x)dsdx\nonumber\\
&=&
\int _{-\infty}^{+\infty}L  _t(x){\rm d}_x\nabla f_n(t,x)-\int _{-\infty}^{+\infty}\int _0^{t}
L  _s(x){\rm d}_{s,x}\nabla f_n(s,x).
\end{eqnarray}
First, consider the case when $\nabla ^-f$ is nondecreasing and left continuous. As we have seen,  $\nabla ^-f$ generates a measure. From
(\ref{zhao33}) it is easy to know that $\nabla f_n$ is also nondecreasing, so also
generates a measure, denoted by $\mu _n$. It is easy to see for any Borel set $G$ in
 $[0,t]\times R^1$,
 \begin{eqnarray*}
 \mu _n(G)&=&\int _{-\infty}^{+\infty}\int _0^t1_{G}{\rm d}_{s,x}\nabla f_n(s,x)
 \\
 &=&\int _0^2\int _0^2\rho(\tau)\rho(z)\int _{-\infty}^{+\infty}\int _0^t1_{G}{\rm d}_{s,x}\nabla ^-f(s-{\tau\over n},x-{z\over n})d\tau dz.
 \end{eqnarray*}
 As $L_s(x)$ is a measurable function so by the definition of Lebesgue
 integrals and substitution of variables, we have
  \begin{eqnarray}\label{convergence1}
&&\int _{-\infty}^{+\infty}\int _0^tL_s(x){\rm d}_{s,x}\nabla f_n(s,x)\nonumber
 \\
 &=&\int _0^2\int _0^2\rho(\tau)\rho(z)\int _{-\infty}^{+\infty}\int _0^tL_s(x){\rm d}_{s,x}\nabla ^-f(s-{\tau\over n},x-{z\over n})d\tau dz\nonumber\\
&=&\int _0^2\int _0^2\rho(\tau)\rho(z)\int _{-\infty}^{+\infty}\int _{-{\tau\over n}}^{t-{\tau\over n}}
L_{s+{\tau\over n}}(x+{z\over n})
{\rm d}_{s,x}\nabla ^-f(s,x)d\tau dz\nonumber\\
&\to & \int _{-\infty}^{+\infty}\int _{0}^{t}
L_{s}(x){\rm d}_{s,x}\nabla ^-f(s,x)\ \ {\rm as} \ \ n\to \infty.
 \end{eqnarray}
 Here the convergence follows by using Lebesgue's dominated convergence theorem
 and noting the local time is right continuous and bounded with compact support almost surely.
 Now for the case when $\nabla ^-f$ is of bounded variation, the above argument applies to us
 by decomposing  $\nabla ^-f$ into the difference of nondecreasing and left continuous functions.
Similarly one can prove that
\begin{eqnarray}\label{convergence2}
\int _{-\infty}^{+\infty}L  _t(x)
{\rm d}_{x}\nabla f_n(t,x)\to
\int _{-\infty}^{+\infty}L  _t(x)
{\rm d}_{x}\nabla^- f(t,x)\ \ {\rm as} \ \ a.s.,
\end{eqnarray}
as $n\to \infty$.
This proves the desired formula. It is noted that (\ref{convergence1}) and (\ref{convergence2})
are also true if we only consider the jump part of the local time,  $\bar L_t(x)$.

To assert $f(t,X(t))$ is a semi-martingale, 
we only need to prove that $\int _{-\infty}^{+\infty}L  _t(x)
{\rm d}_{x}\nabla^- f(t,x)-\int _{0}^{t}\int _{-\infty}^{+\infty}
L_{s}(x){\rm d}_{s,x}\nabla ^-f(s,x)$ is of bounded variation in $t$ on any bounded interval
$[0,T]$.
First we consider the
jump part of the local time. As $\bar
L_t(x)$ is of bounded variation and $\nabla f_n(s,x)$ is smooth, so
$\int _0^t\int _{-\infty}^{\infty} \nabla f_n(s,x){\rm d}_{s,x}\bar
L_{s}(x)$ is a Rieman-Stieltjes integral and the integration by
parts formula holds
\begin{eqnarray}
&&
\int _{-\infty}^{\infty} \bar L_{t}(x){\rm d}_{x}\nabla f_n(t,x)
-\int _0^t\int _{-\infty}^{\infty} \bar L_{s}(x){\rm d}_{s,x}\nabla f_n(s,x)\nonumber\\
&=& -\int _0^t\int _{-\infty}^{\infty} \nabla f_n(s,x){\rm
d}_{s,x}\bar L_{s}(x).
\end{eqnarray}
Then applying Lebesgue's dominated convergence theorem to the latter
integral and (\ref{convergence1}), (\ref{convergence2}) for $\bar L$, we have
\begin{eqnarray}
&&
\int _{-\infty}^{\infty} \bar L_{t}(x){\rm d}_{x}\nabla^- f(t,x)
-\int _0^t\int _{-\infty}^{\infty} \bar L_{s}(x){\rm d}_{s,x}
\nabla^- f(s,x)\nonumber\\
&=& -\int _0^t\int _{-\infty}^{\infty} \nabla^- f(s,x){\rm
d}_{s,x}\bar L_{s}(x).
\end{eqnarray}
It is obvious that $\int _0^t\int _{-\infty}^{\infty} \nabla^-
f(s,x){\rm d}_{s,x}\bar L_{s}(x)$ is of bounded variation in $t$.

Secondly we consider the continuous part of the local time $\tilde L_t(x)$. Define
\begin{eqnarray}
A(t)=\int _{-\infty}^{\infty} \tilde L_{t}(x){\rm d}_{x}\nabla^- f(t,x).
\end{eqnarray}
Assume that $\nabla^- f(t,x)$ is increasing in $x$ at the moment.
For any partition $E_t=\{0=t_0<t_1<\cdots <t_m=T\}$, note first
\begin{eqnarray}\label{convergence3}
A(t_{i+1})-A(t_{i})&=&\int _{-\infty}^{\infty} (\tilde L_{t_{i+1}}(x)-\tilde L_{t_{i}}(x)){\rm d}_{x}\nabla^- f(t_{i+1},x)\nonumber
\\
&&+\int _{-\infty}^{\infty} \tilde L_{t_{i}}(x){\rm d}_{x}(\nabla^- f(t_{i+1},x)-\nabla^- f(t_{i},x))\nonumber\\
:&=&I^1_{t_{i},t_{i+1}}+I^2_{t_{i},t_{i+1}}.
\end{eqnarray}
But  $ I^1_{t_{i},t_{i+1}}$ is a
Rieman-Stieltjes integral, so
\begin{eqnarray}
I^1_{t_{i},t_{i+1}}&=&\lim\limits_{\delta _x\to 0}\sum \limits _{j=0}^{n-1}(\tilde L_{t_{i+1}}(x_j)-\tilde L_{t_{i}}(x_j))
(\nabla^- f(t_{i+1},x_{j+1})-\nabla^- f(t_{i+1},x_{j}))\nonumber\\
&=& \lim\limits_{\delta _x\to 0}\sum \limits _{j=0}^{n-1}(L_{t_{i+1}}(x_j)-L_{t_{i}}(x_j))
(\nabla^- f(t_{i+1},x_{j+1})-\nabla^- f(t_{i+1},x_{j}))\nonumber\\
&&-\lim\limits_{\delta _x\to 0}\sum \limits _{j=0}^{n-1}(\bar L_{t_{i+1}}(x_j)-\bar L_{t_{i}}(x_j))
(\nabla^- f(t_{i+1},x_{j+1})-\nabla^- f(t_{i+1},x_{j}))\nonumber\\
:&=& I^{1a}_{t_{i},t_{i+1}}+I^{1b}_{t_{i},t_{i+1}},
\end{eqnarray}
where $-N=x_0<x_1<\cdots <x_n=N$ is a partition of $[-N,N]$ and
$\delta _x=\max_{j}|x_{j+1}-x_{j}|$ (suppose $N$ is sufficiently
large such that (-N,N) supports $L^{\cdot}_t$). 
We will use the following notation: for any function $g$ of two variables
\begin{eqnarray*} \delta_{ij}g= (g(t_{i+1},x_{j+1})-g(t_{i+1},x_{j})-g(t_{i},x_{j+1})+g(t_{i},x_{j})). \nonumber
\end{eqnarray*}
Noticing
$I^{1a}_{t_{i},t_{i+1}}\geq 0$, and rearranging terms in $
\sum \limits_{i=0}^{m-1}I^{1a}_{t_{i},t_{i+1}}$, we have
\begin{eqnarray*}
\sum \limits_{i=0}^{m-1}I^{1a}_{t_{i},t_{i+1}} =\int
_{-\infty}^{\infty} L_t(x){\rm d}_{x}\nabla^-
f(t,x)-\lim\limits_{\delta _x\to 0}\sum \limits_{i=0}^{m-1} \sum
\limits _{j=0}^{n-1} L_{t_{i}}(x_j)\delta_{ij}\nabla ^-f,
\end{eqnarray*}
As $\nabla ^-f(s,x)$ is of bounded variation in $(s,x)$, thus
\begin{eqnarray*}
\sup_{E_t}\sum \limits_{i=0}^{m-1}I^{1a}_{t_{i},t_{i+1}} &\leq &\int
_{-\infty}^{\infty} L_t(x){\rm d}_{x}\nabla^-
f(t,x)\nonumber\\
&&
+(\sup_{x\in [-N,N]}L_T(x)) \sup_{E_t}\lim\limits_{\delta _x\to
0}\sum \limits_{i=0}^{m-1} \sum \limits _{j=0}^{n-1}
|\delta_{ij}\nabla ^-f|<\infty.
\end{eqnarray*}
Similarly, as $\bar L$ is of bounded variation in $(s,x)$, so
\begin{eqnarray*}
\sup_{E_t}\sum \limits_{i=0}^{m-1}|I^{1b}_{t_{i},t_{i+1}}| &= &
\sup_{E_t}\lim\limits_{\delta _x\to 0}\sum \limits_{i=0}^{m-1} |\sum
\limits _{j=0}^{n-1} \nabla ^-f(t_{i+1},x_{j+1})\delta_{ij}\bar L|\\
&\leq &
(\sup_{0\le s\le T,x\in [-N,N]}|\nabla ^-f(s,x)|) \sup_{E_t}\lim\limits_{\delta _x\to
0}\sum \limits_{i=0}^{m-1} \sum \limits _{j=0}^{n-1}
|\delta_{ij}\bar L|<\infty.
\end{eqnarray*}
 For $ I^2_{t_{i},t_{i+1}}$, noticing $ I^2_{t_{i},t_{i+1}}$ is 
 a Rieman-Stieltjes integral, we also have
\begin{eqnarray*}
\sup_{E_t}\sum \limits_{i=1}^{m-1}|I^2_{t_{i},t_{i+1}}|
=(\sup_{x\in [-N,N]}\tilde L_T(x)) \sup_{E_t}\lim\limits_{\delta _x\to
0}\sum \limits_{i=0}^{m-1} \sum \limits _{j=0}^{n-1}
|\delta_{ij}\nabla ^-f|
<&\infty.
\end{eqnarray*}
This proves that $A(t)$ is of bounded variation. The general case when
$\nabla ^-f(t,x)$ is of bounded variation in $x$ can be proved by
applying the above result for the case to the difference of two increasing functions. The term
$B(t)=\int _0^t\int _{-\infty}^{\infty} \tilde L_{s}(x){\rm d}_{s,x}
\nabla^- f(s,x)$ is obviously of bounded variation, so is
$A(t)+B(t)$. This asserts the claim.
 $\hfill\diamond$
\bigskip

 The smoothing procedure can easily be modified to prove that  if $f:
 R^+\times R\to R$ satisfies (i), (ii) and (iii) of Theorem \ref{tom100},
 is also  $C^1$ in $x$ and the left derivative
 $\Delta ^-f(t,x)$ exists at all points of $[0,\infty)\times R$ and is jointly left continuous and locally
 bounded, then $\Delta f_n(t,x)\to \Delta ^-f(t,x)$ as
 $n\to \infty$, $t>0$. Thus
  \begin{eqnarray}\label{krylov7}
 f(t,X(t))
&=&f(0,X(0))+\int _0^t{\partial ^-\over \partial s} f(s,X(s))ds+\int _0^t\nabla f(s,X(s))dX(s)\nonumber\\
&&
 +{1\over 2}\int _0^t\Delta ^-f(s,X(s))d<X>_s
\ \ a.s.
 \end{eqnarray}

 The next theorem is an easy extension of Theorem \ref{tom100} and formula (\ref{krylov7}).

\begin{thm} \label{tom102}
Assume $f: R^+\times R\to R$ satisfies conditions (i), (ii) and (iii) of Theorem \ref{tom100}.
Further suppose $f(t,x)=f_h(t,x)+f_v(t,x)$ where
\vskip2pt

(i)  $f_h(t,x)$ is $C^1$ in $x$ with
 $\nabla f_h(t,x)$ having left partial derivative
 $\Delta ^-f_h(t,x)$, (with respect to $x$), which is left continuous and locally bounded,
\vskip2pt

 (ii)
 $f_v(t,x)$ has a left continuous derivative $\nabla ^-f_v(t,x)$ at all points $(t,x)$ $[0,\infty)\times R$, which is of locally bounded variation in $(t, x)$ and of locally bounded in $x$ for $t=0$.
\vskip2pt

 \noindent
Then for any continuous semi-martingale $\{X(t), t\geq 0\}$,
$f(t,X(t))$ is a semi-martingale and
\begin{eqnarray}\label{krylov6}
f(t,X(t))
&=&f(0,X(0))+\int _0^t{\partial ^- \over \partial s} f(s,X(s))ds+\int _0^t\nabla^- f(s,X(s))dX(s)\nonumber\\
&&
 +{1\over 2}\int _0^t\Delta ^-f_h(s,X(s))d<X>_s+
\int _{-\infty}^{\infty }L  _t(x){\rm d}_x\nabla  ^-f_v(t,x)\nonumber\\
 &&
-\int _{-\infty}^{+\infty}\int _0^{t}L  _s(x)
{\bf \rm d}_{s,x}\nabla ^- f_v(s,x)\ \ a.s.
\end{eqnarray}
\end{thm}
 \vskip5pt
 {\em Proof}. Mollify $f_h$ and $f_v$, and so $f$, as in the proof of Theorem \ref{tom100}.
 Apply It${\hat {\rm o}}$'s formula to the mollification of $f$ and take the limits as in the proofs of
 Theorem \ref{tom100} and (\ref{krylov7}). \hfill $\diamond$
 \vskip5pt

 If $f$ has discontinuity of first and second order derivatives across a curve $x=l(t)$, where
 $l(t)$ is a continuous function of locally bounded
 variation, it will be convenient to consider the continuous semi-martingale
\begin{eqnarray*}
X^*(s)=X(s)-l(s),
\end{eqnarray*}
 and let $L  _s^*(a)$ be its local time. We can prove the following version of our main results:

\begin{thm}\label{ito2}
 Assume $f: R^+\times R\to R$ satisfies conditions (i), (ii) and (iii) of Theorem \ref{tom100}.
Moreover,  suppose $f(t,x)=f_h(t,x)+f_v(t,x)$, where $f_h(t,x)$ is $C^1$ in $x$ and
 $\nabla f_h(t,x)$ has left derivative
 $\Delta ^-f_h(t,x)$ which is left continuous and locally bounded, and
there exists a curve $x=l(t)$, $t\geq 0$, a continuous
function of locally bounded variation such that
$\nabla  ^-f_v(t,x+l(t))$ as a function of $(t,x)$
is of locally bounded variation in $(t,x)$ and of locally bounded in $x$ for $t=0$. Then
\begin{eqnarray}\label{ito13}
&&
f(t,X(t))\nonumber\\
&=&f(0,z)+\int _0^t{\partial  \over \partial s} f(s,X(s))ds+\int _0^t\nabla^- f(s,X(s))dX(s)\nonumber \\
&&
 +{1\over 2}\int _0^t\Delta ^-f_h(s,X(s))d<X>_s
+
\int _{-\infty}^{\infty }L  _t^*(x){\rm d}_x\nabla  ^-f_v(t,x+l(t))\nonumber\\
 &&
 -\int _{-\infty}^{+\infty}\int _0^{t}L  _s^*(x)
{\bf \rm d}_{s,x}\nabla  ^-f_v(s,x+l(s)))\ \ a.s.
\end{eqnarray}
\end{thm}
\vskip5pt
{\em Proof}. We only need to consider the case
when $f_h=0$ as the general case will follow
easily. We basically follow the proof of Theorem \ref{tom100} and apply It${\hat {\rm o}}$'s
formula to $f_n$ and $X(s)$. We still have (\ref{zhao11}). But by the occupation times formula, a.s.
\begin{eqnarray*}
&&
{1\over 2}\int _0^t\Delta f_n(s,X(s))d<M,M>_s\nonumber\\
&=&{1\over 2}\int _0^t\Delta f_n(s,X^*(s)+l(s))d<M,M>_s\nonumber\\
 &=&\int _{-\infty}^{+\infty}\int _0^t\Delta f_n(s,x+l(s)){\rm d}_sL  ^*_s(x)dx\nonumber\\
&=&
\int _{-\infty}^{+\infty}\Delta f_n(t,x+l(t))L  ^* _t(x)dx-\int _{-\infty}^{+\infty}\int _0^{t}{d \over ds}
\Delta f_n(s,x+l(s))L  ^*_s(x)dsdx\nonumber\\
&\to &
\int _{-\infty}^{\infty }L  ^* _t(x){\rm d}_x\nabla^- f(t,x+l(t))-\int _{-\infty}^{+\infty}\int _0^{t}L  ^*_s(x)
{\bf \rm d}_{s,x}\nabla^- f(s,x+l(s)),
\end{eqnarray*}
as $n\to \infty$ as in the proof of Theorem \ref{tom100}.
This proves the desired formula. \hfill $\diamond$

\begin{cor}
Assume $f: R^+\times R\to R$ satisfies condition (i) of Theorem \ref{tom100} and its left derivative
${\partial ^-\over \partial t}f$ exists on $(0,\infty)\times R$ and is left continuous. Further suppose that there exists a curve $x=l(t)$ of locally
bounded variation such that
f is $C^1$ in x off the curve with $\nabla f$ having left and right limits in $x$ at each point $(t,x)$ and  a left continuous and locally bounded left derivative $\Delta^- f$
on $x$ not equal to $l(t)$. Also assume $\nabla f(t,l(t)+y-)$ as a function of $t$ and $y$
is locally bounded and jointly left continuous if $y\leq 0$, and $\nabla f(t,l(t)+y+)$ is locally
bounded and jointly left continuous in $t$ and right continuous in $y$ if $y\geq 0$.
 Then for any continuous semi-martingale $\{X(t), t\geq 0\}$,
\begin{eqnarray}\label{ito11}
&&
f(t,X(t))\nonumber\\
&=&f(0,X(0))+\int _0^t{\partial  ^-\over \partial s} f(s,X(s))ds+\int _0^t\nabla^- f(s,X(s))dX(s)\nonumber\\
&&
+{1\over 2}\int _0^t\Delta  ^-f(s,X(s))d<X,X>_s\nonumber\\
&&+
\int _0^t(\nabla f(s,l(s)+)-\nabla f(s,l(s)-)){\rm d}_sL  _s^*(0)\ \ a.s.
\end{eqnarray}
\end{cor}
\vskip5pt
{\em Proof}. At first we assume temporarily that $(\nabla f(t,l(t)+)-\nabla f(t,l(t)-))$ is of bounded variation.
This condition will be dropped later.
Formula (\ref{ito11}) can be read from (\ref{ito13}) by considering
 \begin{eqnarray*}
f_h(t,x)&=&f(t,x)+(\nabla f(t,l(t)-)-\nabla f(t,l(t)+))(x-l(t))^+,
\\
f_v(t,x)&=&(\nabla f(t,l(t)+)-\nabla f(t,l(t)-))(x-l(t))^+,
 \end{eqnarray*}
 and integration by parts formula and noticing $\nabla ^-f_v(t,x+l(t))$ is of locally bounded
 variation in $(t,x)$.  Let
$g(t,y)=f(t,y+l(t))$.
 In terms of $X^*$, (\ref{ito11}) can be rewritten as
 \begin{eqnarray}\label{krylov10}
&&
g(t,X^*(t))\nonumber\\
&=&g(0,X^*(0))+\int _0^t g(ds,X^*(s))
+\int _0^t\nabla^- g(s,X^*(s))dX^*(s)\nonumber\\
&&
+{1\over 2}\int _0^t\Delta  ^-g(s,X^*(s))d<X^*,X^*>_s\nonumber\\
&&+
\int _0^t(\nabla g(s,0+)-\nabla g(s,0-)){\rm d}_sL  _s^*(0)\ \ a.s.
\end{eqnarray}
Here
\begin{eqnarray*}
g(ds,y)={\rm d}_sg(s,y)={\partial^-\over \partial s}f(s,y+l(s))ds+\nabla ^-f(s,y+l(s))dl(s).
\end{eqnarray*}

 Now without assuming that $(\nabla f(t,l(t)+)-\nabla f(t,l(t)-))$ is of bounded variation,
 we can prove the formula by a smoothing procedure in the variable $t$. To see this, let
 \begin{eqnarray*}
g_n(t,y)=\int _0^2\rho (\tau)g(t-{\tau\over n},y)d\tau
=\int _0^2\rho (\tau)f(t-{\tau\over n},y+l(t-{\tau\over n}))d\tau,
 \end{eqnarray*}
 with $l(s)=l(0)$ if $s<0$ and $f(s,x)=f(-s,x)$ for $s<0$ as usual.
 Then as $n\to \infty$,
  \begin{eqnarray}\label{ky1}
\int _0^t g_n(ds,X^*(s))
&=& \int _0^t \int _0^2\rho (\tau){\partial^- \over \partial s}f(s-{\tau\over n},X^*(s)+l(s-{\tau\over n}))d\tau ds\nonumber\\
&&+\int _0^t\int _0^2\rho (\tau)\nabla ^-f(s-{\tau\over n},X^*(s)+l(s-{\tau\over n}))dl(s-{\tau\over n})d\tau
\nonumber\\
&\to & \int _0^t{\partial^- \over \partial s}f(s,X^*(s)+l(s))ds
+\nabla ^-f(s,X^*(s)+l(s))dl(s)\nonumber\\
&=&\int _0^tg(ds,X^*(s))\ \ a.s.
 \end{eqnarray}
 It is easy to see that for all $(t,y)$
 \begin{eqnarray}\label{ky2}
 g_n(t,y)\to g(t,y)
 \end{eqnarray}
 and
 for all $y\ne 0$,
 \begin{eqnarray}\label{ky3}
 \nabla g_n(t,y)\to \nabla g(t,y),  \Delta ^-g_n(t,y)\to \Delta ^-g(t,y),
 \end{eqnarray}
 with uniform local bounds.
 Moreover, we can see that as  $y\to 0\pm$ and $n\to \infty$,
  \begin{eqnarray}\label{ky4}
\nabla ^{\pm}g_n(t,y)=\int _0^2\rho (\tau)\nabla ^{\pm}g(t-{\tau\over n},y)d\tau \to \nabla g(t,0\pm).
\end{eqnarray}

 Since $\nabla g_n(t,0\pm )$ are smooth in $t$ then are of locally bounded variation.
 From (\ref{krylov10}),
  \begin{eqnarray}\label{krylov11}
&&
g_n(t,X^*(t))\nonumber\\
&=&g_n(0,X^*(0))+\int _0^t g_n(ds,X^*(s))
+\int _0^t\nabla^- g_n(s,X^*(s))dX^*(s)\nonumber\\
&&
+{1\over 2}\int _0^t\Delta  ^-g_n(s,X^*(s))d<X^*,X^*>_s\nonumber\\
&&+
\int _0^t(\nabla g_n(s,0+)-\nabla g_n(s,0-)){\rm d}_sL  _s^*(0)\ \ a.s.
\end{eqnarray}
We obtain the desired formula by passing to
the limits using (\ref{ky1}), (\ref{ky2}), (\ref{ky3}) and (\ref{ky4}).
 $\hfill\diamond$

 \begin{rmk} (i) Formula (\ref{ito11}) was also observed by Peskir in \cite{pe}
 and \cite{gp} independently.

(ii) From the proof of Theorem \ref{tom100}, one can take different mollifications,
e.g. one can take (\ref{truman41}) as
\begin{eqnarray*}
f_n(s,x)=\int _0^2\int _0^2\rho(\tau)\rho(z)f(s+{\tau\over n},x+{z\over n})d\tau dz, \ \ n\geq 1.
\end{eqnarray*}
This will lead to as $n\to \infty$,
\begin{eqnarray*}
{\partial\over \partial s} f_n(s,x)\to {\partial^+\over \partial s} f(s,x)
\end{eqnarray*}
instead of (\ref{zhao2}), if $ {\partial^+\over \partial s} f(s,x)$ is jointly right continuous.
 Therefore we have the following more general It$\hat  o$'s formula
\begin{eqnarray*}\label{truman42}
f(t,X(t))
&=&f(0,z)+\int _0^t{\partial ^{s_1}\over \partial s} f(s,X(s))ds+\int _0^t\nabla ^{s_2}f(s,X(s))dX(s)\nonumber\\
&&
+{1\over 2}\int _0^t\Delta ^{s_2} f_h(s,X(s))d<X>_s\nonumber\\
&&
+\int _{-\infty}^{\infty }L  _t^{s_2}(x){\rm d}_x\nabla ^{s_2}f_v(t,x)
-\int _{-\infty}^{+\infty}\int _0^{t}L  _s^{s_2}(x)
{\bf \rm d}_{s,x}\nabla ^{s_2}f_v(s,x)\ \ a.s.,
\end{eqnarray*}
where $s_1=\pm$ and $s_2=\pm $, $L  _s^{s_2}(x)$ is taken to be the cadlag version when $s_2=-$ and the caglad version
otherwise.
\end{rmk}
\vskip10pt

 Formula (\ref{krylov6}) is in a very general form. It includes
the  classical It$\hat{\rm o}$ formula, Tanaka's formula,
 Meyer's formula for convex functions, the formula given by Az\'ema, Jeulin, Knight and Yor \cite{azema}
 and formula (\ref{ito11}).   In the following we will give some
 examples for which (\ref{ito11}) and some known  generalized It${\hat {\rm o}}$ formulae do
 not immediately apply, but formula (\ref{krylov6}) can be applied.
 These examples can be presented in different forms to include local times on curves.
 \vskip10pt

 \begin{exap}\label{example1} Consider the function
 \begin{eqnarray*}
 f(t,x)=(\sin \pi x\sin\pi t)^+.
 \end{eqnarray*}
 Then
 \begin{eqnarray*}
 \nabla ^-f(t,x)=\pi \cos \pi x\sin\pi t1_{\sin \pi x\sin\pi t>0}.
 \end{eqnarray*}
 One can verify that $\nabla ^-f(t,x)$ is of locally bounded variation in $(t,x)$. This can be easily seen
 from Proposition \ref{prop1} and the simple fact that
\begin{eqnarray*}
&&  \cos \pi x\sin\pi t 1_{\sin \pi x\sin \pi t>0}\\
& =&
\left\{ \begin{array}{l}
\cos \pi x\sin\pi t, \ {\rm if}\  i\leq t< i+1,
 j\leq x< j+1, i+j \ {\rm is \ even}
\\
0, \ \ \ \ \ \ \ \ \ \ \ \ \ \ \
{\rm otherwise}
 \end{array} \right.
\end{eqnarray*}
 Therefore
 \begin{eqnarray*}
 (\sin \pi X(t)\sin\pi t)^+&=&\pi \int _0^t \cos \pi s\sin\pi X(s)1_{\sin \pi X(s)\sin\pi s>0}ds
 \\
 &&+\pi \int _0^t \cos \pi X(s)\sin\pi s1_{\sin \pi X(s)\sin\pi s>0}dX(s)\\
 &&
+\pi \sin \pi t \int _{-\infty}^{\infty}L  _t(a){\rm {\rm d}_a}(\cos \pi a 1_{\sin \pi a\sin\pi t>0})\\
&&-\pi
\int _0^t \int _{-\infty}^{\infty}L  _s(a){\rm {\rm d}_{s,a}}(\cos \pi a\sin\pi s1_{\sin \pi a\sin\pi s>0}).
 \end{eqnarray*}
 One can expand the last two integrals to see the jump of $\cos \pi a\sin\pi s1_{\sin \pi a\sin\pi s>0}$.
  \end{exap}

 Note in example \ref{example1}, $\nabla ^-f(t,x)$ has jump on the boundary of each interval
 $i\leq t<i+1$, $j\leq x< j+1$. One can use this example as a prototype to construct
 many other examples with other types of derivative jumps.

 \begin{exap} Consider the function
 \begin{eqnarray*}
 f(t,x)=(\sin \pi x)^{1\over 3} (\sin \pi x\sin\pi t)^+.
 \end{eqnarray*}
 Then
 \begin{eqnarray*}
 \nabla ^-f(t,x)&=&{1\over 3}\pi \cos \pi x(\sin\pi x)^{-{2\over 3}} (\sin \pi x\sin\pi t)^+\\
 && +
 \pi (\sin \pi x)^{1\over 3}\cos \pi x\sin\pi t1_{\sin \pi x\sin\pi t>0}.
 \end{eqnarray*}
 One can verify that $\nabla ^-f(t,x)$ is of locally bounded variation in $(t,x)$ and continuous. In fact,
 \begin{eqnarray*}
&&
 \cos \pi x(\sin\pi x)^{-{2\over 3}} (\sin \pi x\sin\pi t)^+ \\
 &=&\left\{ \begin{array}{l}
 \cos \pi x(\sin\pi x)^{{1\over 3}} \sin\pi t , \ {\rm if}\  i\leq t< i+1,\\
 \ \ \ \ \ \ \ \ \ \ \ \ \ \ \ \ \ \ \ \ \ \ \ \ \ \ \ \ j\leq x< j+1, i+j \ {\rm is \ even}
\\
0, \ \ \ \ \ \ \ \ \ \ \ \ \ \ \ \ \ \ \ \ \ \ \ \ \
{\rm otherwise}.
 \end{array} \right.
\end{eqnarray*}
then it is easy to see that $\cos \pi x(\sin\pi x)^{-{2\over 3}} (\sin \pi x\sin\pi t)^+$ is
of locally bounded variation in $(t,x)$ using proposition \ref{prop1}. Similarly one can see that
$(\sin \pi x)^{1\over 3}\linebreak \cos \pi x \sin\pi t1_{\sin \pi x\sin\pi t>0}$ is of locally
bounded variation in $(t,x)$ as well.

 Note
 $\Delta ^-f(t,x)$ blows up when $x$ is near an integer value, and their left and right limits
 also blow up. However one can apply our generalized It${\hat{\rm o}}$'s formula (\ref{krylov6})
 to this function so that
 \begin{eqnarray*}
 &&
 (\sin \pi X(t))^{1\over 3} (\sin \pi X(t)sin\pi t)^+\nonumber\\
 &=&\pi \int _0^t  (\sin \pi X(s))^{4\over 3} \cos \pi s 1_{\sin \pi X(s)\sin\pi s>0}{\rm ds}\nonumber\\
&& +\int_{-\infty}^{\infty}L_t(a){\rm {\rm d}_a}({1\over 3}\pi \cos \pi a(\sin\pi a)^{-{2\over 3}} (\sin \pi a\sin\pi t)^+\nonumber \\
&& \hskip3cm +
 \pi (\sin \pi a)^{1\over 3}\cos \pi a\sin\pi t1_{\sin \pi a\sin\pi t>0})\nonumber\\
 &&
 -\int _0^t\int _{-\infty}^{\infty}L_s(a){\rm {\rm d}_{s,a}}({1\over 3}\pi \cos \pi a(\sin\pi a)^{-{2\over 3}} (\sin \pi a\sin\pi s)^+\nonumber\\
 &&\hskip3cm +
 \pi
 (\sin \pi a)^{1\over 3} \cos \pi a\sin\pi s1_{\sin \pi a\sin\pi s>0}).
 \end{eqnarray*}

  \end{exap}

\section{The case for It${\hat {\rm o}}$ processes}

\setcounter{equation}{0}

For It${\hat {\rm o}}$ processes, we can allow some of the generalized derivatives of $f$ to be
only in
$L^2_{loc}(dtdx)$. Consider
\begin{eqnarray}\label{krylov8}
X(t)=X(0)+\int_0^t\sigma_rdW_r+\int _0^tb_rdr.
\end{eqnarray}
Here $W_r$ is a one-dimensional Brownian motion on a filtered probability space $(\Omega,{\cal F},
\{{\cal F}_r\}_{r\geq 0}, P)$ and $\sigma_r$ and $b_r$ are progressively measurable with respect to
$\{{\cal F}_r\}$ and satisfy the following conditions:
for all $t>0$
\begin{eqnarray}\label{krylov1}
\int _0^t|\sigma _r|^2dr<\infty, \ \ \int _0^t|b _r|dr<\infty \ \ a.s.
\end{eqnarray}
Under condition (\ref{krylov1}), the process (\ref{krylov8}) is well defined.
For any $N>0$, define $\tau_N=\inf\{s:|X(s)|\geq N\}$. Assume
there exist constants $\delta >0$ and $K>0$
such that  ,
\begin{eqnarray}\label{krylov9}
\sigma _t(\omega)\geq \delta >0,\ \ |\sigma _t(\omega)|+|b_t(\omega)|\leq K,\ \
{\rm for\ all}\  (t,\omega) {\rm \ with \ } t\leq \tau _N.
\end{eqnarray}
The following inequality due to Krylov \cite{krylov} plays an important role.

\begin{lem}\label{krylov2}
Assume condition (\ref{krylov1}) and (\ref{krylov9}). Then there exists a constant $M>0$, depending
only on $\delta$ and $K$ such that
\begin{eqnarray}\label{krylov4}
E\int _0^{t\wedge\tau_N} |f(r,X(r))|dr\leq M(\int_0^t\int _{-N}^{+N}(f(r,x))^2drdx)^{1\over 2}.
\end{eqnarray}
\end{lem}

Denote again by $L_t(x)$ the local time of the diffusion process $X(t)$ at level $x$.
We can prove the following theorem.

\begin{thm}\label{krylov8}
Assume $f(t,x)$ is continuous with
generalized derivative ${\partial \over \partial t}f$  in $L^2_{loc}(dtdx)$
and generalized derivative $\nabla f$ of
 locally bounded variation in $(t, x)$ and of locally bounded variation in $x$
 for $t=0$. Consider  an It${\hat {\rm o}}$
 process $X(t)$ given by (\ref{krylov8}) with $\sigma$ and $b$ satisfying (\ref{krylov1}) and (\ref{krylov9}).
Then  a.s.
\begin{eqnarray}
f(t,X(t))
&=&f(0,X(0))+\int _0^t{\partial \over \partial s} f(s,X(s))ds+\int _0^t
 \nabla f(s,X(s))dX(s)\nonumber\\
&&
+
\int _{-\infty}^{\infty }L  _t(x){\rm d}_x\nabla f(t,x)-\int _{-\infty}^{+\infty}\int _0^{t}L  _s(x)
{\bf \rm d}_{s,x}\nabla f(s,x).
\end{eqnarray}
\end{thm}
 \vskip5pt
{\em Proof}. Define $f_n$ by (\ref{krylov5}). From a well-known result on Sobolev
spaces (see Theorem 3.16, p.52 in \cite{adams}),  we know that as  $n\to \infty$,
\begin{eqnarray*}
f_n(t,x) &\to & f(t,x) ,
\end{eqnarray*}
for all $(t,x)$ and for any $N>0$
\begin{eqnarray*}
{\partial \over \partial t}f_n&\to & {\partial  \over \partial t} f,\label{zhao02} {\rm \ in}\  L^2([0,t]\times [-N,N])\\
\nabla f_n&\to & \nabla f,\label{zhao03} {\rm \ in}\   L^4([0,t]\times [-N,N]).
\end{eqnarray*}

As in the proof of Theorem \ref{tom100}, we have the It$\hat {\rm o}$ formula (\ref{zhao11}) for
$f_n(t\wedge \tau _N, X(t\wedge \tau _N))$. The convergence of the terms
$f_n(t\wedge \tau _N, X(t\wedge \tau _N))$, and ${1\over 2}
\int _0^{t\wedge \tau _N}\sigma _s^2\Delta f_n(s,X(s))ds$ is the same as before. Now by using Lemma \ref{krylov2},
\begin{eqnarray*}
&&
E|\int _0^{t\wedge\tau_N}{\partial \over \partial s} f_n(s,X(s))ds-\int _0^{t\wedge\tau_N}{\partial \over \partial s} f(s,X(s))ds|\nonumber\\
&\leq &
E\int _0^{t\wedge\tau_N}|{\partial \over \partial s} f_n(s,X(s))-{\partial \over \partial s} f(s,X(s))|ds\nonumber\\
&\leq &
M(\int _0^t\int _{-N}^{N}({\partial \over \partial s} f_n(s,x)-{\partial \over \partial s} f(s,x))^2dsdx)^{1\over 2}\to 0
\end{eqnarray*}
as $n\to \infty$.
Similarly one can prove
\begin{eqnarray*}
\int _0^{t\wedge\tau_N} b_s\nabla f_n(s,X(s))ds\to \int _0^{t\wedge\tau_N}b_s\nabla f(s,X(s))ds \ \ in \ \ L^1(dP).
\end{eqnarray*}
Moreover, there exists a constant $M>0$ such that
\begin{eqnarray*}
&&
E(\int _0^{t\wedge\tau_N}\sigma _s \nabla f_n(s,X(s))dW_s-\int _0^t\sigma _s\nabla f(s,X(s))dW_s)^2\nonumber
\\
&=&
E(\int _0^{t\wedge\tau_N}\sigma _s ^2(\nabla f_n(s,X(s))-\nabla f(s,X(s)))^2ds\nonumber\\
&\leq &
M(\int _0^t\int _{-N}^{N}(\nabla f_n(s,x)-\nabla f(s,x))^4dsdx\to 0
\end{eqnarray*}
as $n\to \infty$.
Therefore we have proved that
\begin{eqnarray*}
&&
f(t\wedge\tau_N,X(t\wedge\tau_N))\nonumber\\
&=&f(0,X(0))+\int _0^{t\wedge\tau_N}{\partial \over \partial s} f(s,X(s))ds+\int _0^{t\wedge\tau_N}
 \nabla f(s,X(s))dX(s)\\
&&
+
\int _{-\infty}^{\infty }L  _{t\wedge\tau_N}(x){\rm d}_x\nabla f(t\wedge\tau_N,x)-\int _{-\infty}^{+\infty}\int _0^{t\wedge\tau_N}L  _s(x)
{\bf \rm d}_{s,x}\nabla f(s,x).\nonumber
\end{eqnarray*}
The desired formula follows.
$\hfill\diamond$
\bigskip

Recall the following extension of It${\hat {\rm o}}$'s formula due to Krylov (\cite{krylov}):
   if $f:
 R^+\times R$ is $C^1$ in $x$ and $\nabla f$ is absolutely continuous with respect
 to $x$ for each $t$ and the generalized derivatives ${\partial \over \partial s}f(s,x)$ and
 $\Delta f$ are in $L^2_{loc}(dtdx)$, then
 \begin{eqnarray}\label{ito14}
 f(t,X(t))
&=&f(0,z)+\int _0^t{\partial \over \partial s} f(s,X(s))ds+\int _0^t\nabla f(s,X(s))dX(s)\nonumber\\
&&
 +{1\over 2}\int _0^t\sigma _s^2\Delta f(s,X(s))ds
\ \ a.s.
 \end{eqnarray}

 The next theorem is an easy consequence of the method of proof of
 Theorem \ref{krylov8} and of formula (\ref{ito14}).

\begin{thm}
Assume $f(t,x)$ is continuous and
its generalized derivative ${\partial \over \partial t} f$ is in $L^2_{loc}(dtdx)$.
Moreover $f(t,x)=f_h(t,x)+f_v(t,x)$ with $f_h(t,x)$ being $C^1$ in $x$ and
 $\nabla f_h(t,x)$ having generalized derivative
 $\Delta f_h(t,x)$ in $L^2_{loc}(dtdx)$, and $f_v$ having generalized derivative
 $\nabla f_v(t,x)$ being of locally bounded variation in $(t, x)$ and of locally
 bounded variation in $x$ for $t=0$. Suppose $X(t)$ is an It${\hat {\rm o}}$
 process given by (\ref{krylov8})
with $\sigma$ and $b$ satisfying (\ref{krylov1}) and (\ref{krylov9}).
Then,
\begin{eqnarray}\label{ito12}
f(t,X(t))
&=&f(0,X(0))+\int _0^t{\partial \over \partial s} f(s,X(s))ds+\int _0^t\nabla f(s,X(s))dX(s)\nonumber\\
&&
 +{1\over 2}\int _0^t\Delta f_h(s,X(s))d<X>_s+
\int _{-\infty}^{\infty }L  _t(x){\rm d}_x\nabla  f_v(t,x)\nonumber\\
 &&
-\int _{-\infty}^{+\infty}\int _0^{t}L  _s(x)
{\bf \rm d}_{s,x}\nabla f_v(s,x)\ \ a.s.
\end{eqnarray}
\end{thm}

  \bigskip

{\bf Acknowledgment}
\vskip5pt

\footnotesize

It is our great pleasure to thank M. Chen, M. Freidlin, Z. Ma, S. Mohammed, B. \O ksendal, S. Peng,
L.M. Wu, J.A. Yan,
M. Yor and W.A. Zheng for useful conversations.
 The first version of this paper was presented in the Kautokeino
Stochastic Analysis Workshop in July 2001 organized by B. \O ksendal.
We would also like to thank S. Albeverio, Z. Ma,
M. R\"ockner for inviting us to the Sino-German Stochastic Analysis
Meeting (Beijing 2002). HZ would like to thank T.S. Zhang and
T. Lyons for invitations to Manchester and Oxford respectively to
present the results of this paper.  One version of this paper was also presented in
Swansea Workshop on Probabilistic Methods in Fluid in April 2002 and the final  version
in Warwick SPDEs Workshop in August 2003 and Mini-workshop of Local Time-Space
Calculus with Applications in Obverwolfach in May 2004. We would like to thank
G. Peskir and N. Eisenbaum for invitations to the Oberwolfach conference.
We would like to thank Y. Liu, C.R. Feng and B. Zhou for reading the manuscript and
making some valuable suggestions. It is our pleasure to thank the referee
for useful comments. This project is partially supported by EPSRC grants
GR/R69518 and GR/R93582.

\end{document}